\documentclass[11pt]{article}
\usepackage[dvipsnames]{xcolor}
\usepackage{amsfonts}
\usepackage{setspace}
\usepackage{bm}
\usepackage{bbm}
\usepackage{amssymb}
\usepackage{indentfirst}
\usepackage{geometry}
\usepackage{mathtools}
\usepackage{hyperref}
\hypersetup{colorlinks=false, hidelinks, urlcolor=cyan}

\DeclarePairedDelimiter\abs{\lvert}{\rvert}
\DeclarePairedDelimiter{\floor}{\lfloor}{\rfloor}

\DeclareMathOperator{\HP}{HP}

\makeatletter
\let\oldabs\abs
\def\abs{\@ifstar{\oldabs}{\oldabs*}}

\newcommand*{\citena}[1]{%
\begingroup
\color{Green}
\romannumeral-`\x 
\setcitestyle{numbers}%
\citep{#1}
\endgroup
\ignorespacesafterend
}

\newcommand*{\citesup}[1]{%
\begingroup
\color{Green}
\citep{#1}
\endgroup
\ignorespacesafterend
}

\newcommand*{\eqrefe}[1]{%
\begingroup
(\color{BrickRed}
\romannumeral-`\x 
\setcitestyle{numbers}%
\ref{eq:#1}%
\endgroup
)\ignorespacesafterend
}

\newcommand*{\secrefe}[1]{%
\begingroup
(\color{Aquamarine}
\romannumeral-`\x 
\setcitestyle{numbers}%
\ref{#1}%
\endgroup
)\ignorespacesafterend
}

\newtagform{Tags}[\textcolor{BrickRed}]{\color{Black}(}{)}

\usepackage[numbers, sort]{natbib}
\setcitestyle{super}

\newcommand{\ii}{\bm{i}}

\begin{document}
\title{A Generalization of the Riemann Functional Equation with Three Variables}

\author{
José Risomar Sousa\\
\href{https://orcid.org/0009-0002-7779-0808}
     {\textcolor{magenta}{ORCID: 0009-0002-7779-0808}}
}
     
\date{August 20, 2026}
\maketitle

\begin{abstract}
The discovery that the generalized Riemann functional equation from the literature can be extended to real $b>1$, by exploiting symmetries of a new formula for the Hurwitz zeta function with respect to the shift parameter $(b)$, led to the search for a similar relation for the Lerch $\Phi$ function, $\Phi(e^z,k,b)$, where a new formula is used as a starting point and the symmetries with respect to $b$ are exploited again. This search resulted in a generalization of the Riemann functional equation with three variables, whose domain is then extended to $b>1$, with the aid of a new self-contained Lerch $\Phi$ formula, making the relation hold for all real $b$. As a secondary objective, we address the problem of obtaining the analytic continuation of a finite summation that appears in a formula for the Lerch $\Phi$ function. Above all, this paper demonstrates, with clear elementary mathematics, how a functional equation can be derived.
\end{abstract}

\usetagform{Tags}

\tableofcontents


\section{Introduction}

The paper \citena{GHN}, \textit{Generalized Harmonic Numbers}, introduced a new approach that enabled a new formula to be derived for the generalized harmonic numbers, $H_k(n)$. That paper showed a detailed demonstration, based on power series, of why the formula holds and how it changes depending on the choice of ``initial equation" (e.g., the main ones used were $\sin{\pi j}=0$, $\sin{2\pi j}=0$ and $\cos{2\pi j}=1$, with $j$ an integer).\\

None of the various subsequent papers on the theme, that followed, made it clear that the generalization of that approach,
\begin{equation} \label{eq:Init_eq}
e^{z\,j} = \sum_{q=0}^{\infty}\frac{(z\,j)^q}{q!} \text{,}
\end{equation}
\noindent (which holds for all $z$ and $j$, as the exponential function is analytic everywhere), leads to the same formula for the partial sums of the polylogarithm that was derived in paper \citena{Lerch_pos} using a different path, namely,
\begin{multline} \label{eq:Poly_parc}
\sum_{j=1}^{n}\frac{e^{z\,j}}{j^{k}}=\frac{e^{z\,n}}{2n^{k}}-\frac{1}{2n^{k}}\sum_{j=0}^{k}\frac{(z\,n)^j}{j!}+\sum_{j=1}^{k}\frac{z^{k-j}}{(k-j)!}\sum _{q=1}^{n}\frac{1}{q^{j}}\\
+\frac{z^{k}}{2(k-1)!}\int_0^1(1-u)^{k-1}\left(e^{z\,n\,u}-1\right)\coth{\frac{z u}{2}}\,du \text{,}
\end{multline}
\noindent where each case $k$ is obtained from \eqrefe{Init_eq} --- by dividing each side of the equation successively by $j$ and summing it over $j$, a recurrence is obtained (the integral is the part of the power series thus formed that does not admit a closed form).\\

Note that \eqrefe{Poly_parc} can be proved by applying the delta operator with respect to $n$, $\Delta f(n)=f(n)-f(n-1)$, to both sides of the equation, in which case the resulting integral has a closed form\footnote{Obtained via software.}, which confirms the identity. Hence, the equation holds for every $n$, not just integer $n$.\\

The challenge that is now addressed is how to analytically continue the equation \eqrefe{Poly_parc} to complex $k$. The formula as stated breaks down, for example, when $k$ is a non-positive integer as the integral does not converge (that is, that term is no longer given by the anti-derivatives of the $\mathbbm{1}_{k|n}$ function, as per paper \citena{GHN}, but by the derivatives). The analysis of this problem led to the realization that, by exploiting the symmetries of a Lerch $\Phi$ formula that stems from \eqrefe{Poly_parc}, it is possible to produce a generalization of the Riemann functional equation with three parameters. The domain of that equation can then be extended to $b>1$, making it hold for all real $b$. To achieve that, a new self-contained Lerch $\Phi$ formula is carefully derived through analytic continuation of finite summations.\\

In addition, appendix \secrefe{RFE} presents a review and a restatement of a newly found extension of the generalization of the Riemann functional equation from the literature, which holds when the shift parameter ($b$) is a non-negative real number.\\

This manuscript focuses primarily on real domains of the variables $k$, $z$ and $b$ (in the cases where the focus shifts to purely imaginary $z$, for instance, a transformation is applied so that $z$ can still be assumed to be real).\\

Despite references to the author's previous papers, this paper is sufficiently self-contained and does not require extensive prior reading. It presents an exposition of the main results, while leaving other possible derivations for a possible future paper (such as, the functional equation that is obtained via differentiation with respect to the parameter $z$, instead of $b$). A great effort was put into making this manuscript organized and easy to follow, though many readers may still find it a bit challenging.

\subsection{Generalization to harmonic progressions}

The approach described previously can be generalized to produce formulae for the generalized harmonic progressions, a term that has been used previously to refer to the finite sum below,
\begin{equation} \nonumber
\HP_k(n)=\sum_{j=1}^{n}\frac{1}{(a j+b)^k} \text{,}
\end{equation}
\noindent which yields the partial sums of the Hurwitz zeta function.\\

A few different formulae have been created for $\HP_k(n)$ in previous research, by changing the assumptions on the parameters of the initial equation.\\

It is trivial to state a generalization of the result of the previous section (just by making $j \mapsto j+b$ in \eqrefe{Init_eq}). It then gives the formula below for the partial sums of the Lerch $\Phi$ function, 
\begin{multline} \label{eq:Lerch_partial} 
\sum _{j=1}^{n}\frac{e^{z(j+b)}}{(j+b)^{k}}=-\frac{e^{z\,b}}{2\,b^{k}}+\frac{e^{z(n+b)}}{2(n+b)^{k}}+\frac{1}{2\,b^{k}}\sum_{j=0}^{k}\frac{(z\,b)^j}{j!}-\frac{1}{2(n+b)^{k}}\sum_{j=0}^{k}\frac{(z (n+b))^j}{j!}\\+\sum_{j=1}^{k}\frac{z^{k-j}}{(k-j)!}\sum_{q=1}^{n}\frac{1}{(q+b)^j}+\frac{z^{k}}{2(k-1)!}\int_0^1(1-u)^{k-1}\left(e^{z(n+b) u}-e^{z b u}\right)\coth{\frac{z u}{2}}\,du \text{,}
\end{multline}
\noindent which was also derived using a different path in \citena{Lerch_pos}. It holds for every positive integer $k$ and $n$ and every complex $z$ and $b$.

\section{A new approach for the Lerch $\Phi$}

There is a new approach to obtain the limit of \eqrefe{Lerch_partial}. Note that with the generalization from \eqrefe{Poly_parc}, $n$ no longer has to be integer. This allows a new formula for the Lerch $\Phi$ function to be deduced, where functions of $n$ in \eqrefe{Poly_parc} are replaced by their analytic continuations, as below (note $n$ was replaced by $b$ to keep the usual notation):
\begin{multline} \label{eq:Lerch_full_new}
-e^{z\,b}\,\Phi(e^z,k,b)+\mathrm{Li}_{k}(e^z)=-\frac{e^{z\,b}}{2\,b^{k}}-\frac{1}{2\,b^{k}}\sum_{j=0}^{k}\frac{(z\,b)^j}{j!}+\sum_{j=1}^{k}\frac{z^{k-j}}{(k-j)!}H_j(b)\\
+\frac{z^{k}}{2(k-1)!}\int_0^1(1-u)^{k-1}\left(e^{z\,b\,u}-1\right)\coth{\frac{z u}{2}}\,du \text{}
\end{multline}
\indent Like in the original formula (which reappears in the next section), $z$ can be any complex number, except if both $\Re(z) \ge 0$ and $\abs{\Im(z)}>2\pi$. And noting that if $j \neq 1$ then\citesup{Apostol}\citesup{DLMF},
\begin{equation} \label{eq:Hurwitz_zeta_rel}
H_j(b)=\zeta(j)-\zeta(j,b+1) \text{,}
\end{equation}
\noindent it follows that the second summation on the right-hand side of the equation above can be rewritten as:
\begin{multline} \label{eq:sum_Hk(n)}
\sum_{j=1}^{k}\frac{z^{k-j}}{(k-j)!}H_j(b)=\lim\sum_{j=1}^{k}\frac{z^{k-j}}{(k-j)!}\left(\zeta(j)-\zeta(j,b+1)\right)\\=-\frac{z^k}{k!}+\frac{1}{b^{k}}\sum_{j=0}^{k}\frac{(z\,b)^j}{j!}+\frac{z^{k-1}}{(k-1)!}\left(H(b)-\frac{1}{b}\right)-\sum_{j=2}^{k}\frac{z^{k-j}}{(k-j)!}\left(\zeta(j,b)-\zeta(j)\right)
\end{multline}

\subsection{Transforming the original Lerch $\Phi$} \label{Limit_weak}

The limit of the partial Lerch $\Phi$ formula from \eqrefe{Lerch_partial} was first found in \citena{Lerch_pos}, and the new formula in \eqrefe{Lerch_full_new} is evidence of its correctness (an actual proof is given in section \secrefe{Lerch_proof}). The original formula is given by:
\begin{multline} \label{eq:Lerch_full} 
e^{z\,b}\,\Phi(e^z,k,b)=\frac{e^{z\,b}}{2\,b^{k}}+\sum_{j=2}^{k}\frac{z^{k-j}}{(k-j)!}\left(-\frac{1}{2\,b^j}+\zeta(j,b)\right)\\+\frac{\pi\,z^k}{2(k-1)!}\cot{\pi b}-\frac{z^{k-1}}{(k-1)!}\log{\left(-\frac{z}{2\pi}\right)}\\-\frac{z^{k}}{2(k-1)!}\int_0^1(1-u)^{k-1}\,e^{z\,b\,u}\coth{\frac{z u}{2}}+\frac{2\pi}{z}\left(-1+\frac{\sin{2\pi b u}}{\sin{2\pi b}}\right)\cot{\pi u}\,du \text{}
\end{multline}
\indent Paper \citena{Lerch_pos} mentioned that to find the limit of the partial Lerch $\Phi$ from \eqrefe{Lerch_partial} it was necessary to find an approximation for the harmonic progression of order one (which does not converge). That is, for sufficiently large $n$, one can write the interesting asymptotic expression\citesup{Abra}\citesup{DLMF}:
\begin{equation} \label{eq:HP(n)_assimp}
\sum_{j=1}^{n}\frac{1}{j+b} \sim H(n)-\frac{1}{b}\sum_{j=2}^{\infty}(-b)^j\,\zeta(j)=H(n)-H(b) \text{}
\end{equation}
\indent But this approximation is not a new result, it is in the literature. To prove it, the identity below, that relates the ordinary generating function of the zeta function at the odd integers to the harmonic numbers, can be used, 
\begin{multline} \label{eq:soma_zetas_impares}
\sum_{j=1}^{\infty}b^{2j}\,\zeta(2j+1)=-\pi\int_{0}^{1}\left(\frac{\sin{2\pi b u}}{\sin{2\pi b}}-u\right)\cot{\pi u}\,du \\ = \frac{1}{2\,b}-H(b)-\frac{\pi}{2}\cot{\pi b} = -\frac{1}{2\,b}-H(-b)+\frac{\pi}{2}\cot{\pi b} \text{,}
\end{multline}
\noindent where the rightmost formula follows from a property of odd functions, $f(x)=-f(-x)$.\\

With this identity, formula \eqrefe{Lerch_full} can be transformed into the simpler form below, which also holds for integer or half-integer $b$:
\begin{multline} \label{eq:Lerch_full_alt} 
e^{z\,b}\,\Phi(e^z,k,b)=\frac{e^{z\,b}}{2\,b^{k}}+\sum_{j=2}^{k}\frac{z^{k-j}}{(k-j)!}\left(-\frac{1}{2\,b^j}+\zeta(j,b)\right)\\-\frac{z^{k-1}}{(k-1)!}\left(-\frac{1}{2\,b}+H(b)+\log{\left(-\frac{z}{2\pi}\right)}\right)\\-\frac{z^{k}}{2(k-1)!}\int_0^1(1-u)^{k-1}\,e^{z\,b\,u}\coth{\frac{z u}{2}}-\frac{2\pi}{z}(1-u)\cot{\pi u}\,du \text{}
\end{multline}
\indent Note that the $H(b)$ term can be written as an infinite series in two ways,
\begin{equation} \label{eq:alternating_zeta}
H(b) = -\sum_{j=2}^{\infty}(-b)^{j-1}\,\zeta(j) = \sum_{j=1}^{\infty}\left(\frac{1}{j}-\frac{1}{j+b}\right) \text{,}
\end{equation}
\noindent where the second way stems from the relation \eqrefe{Hurwitz_zeta_rel} at the limit, but also from \eqrefe{HP(n)_assimp} (and it explains why $H(b)$ has singularities at the negative integers).

\subsection{A proof of the Lerch $\Phi$ formula} \label{Lerch_proof}

Although not all the steps are provided, the equation \eqrefe{Lerch_full_alt} can be proved by taking the case $k=1$, for which the integral has a known closed form, according to Mathematica\footnote{Some of the results introduced herein were obtained with the aid of software.}: 
\begin{multline} \nonumber
\int_{0}^{1}e^{z\,b(1-u)}\coth{\frac{z u}{2}}-\frac{2\pi\,e^{z\,b}}{z}(1-u)\cot{\pi u}\,du \\ = -\frac{1}{z}\left(\Phi(e^z,1,-b)+e^{z}\,\Phi(e^z,1,-b+1)\right)-\frac{e^{z\,b}}{z}\left(\frac{1}{b}+2\log{\left(-\frac{z}{2\pi}\right)}+2\,H(-b)\right) \text{,}
\end{multline}
\noindent which proves that the formula holds for $k=1$.\\

To prove the general case, the integral above is differentiated $k-1$ times with respect to $b$. The only challenging term is given below:
\begin{equation} \nonumber
\frac{\text{d}^{k-1}(e^{z\,b}\,H(-b))}{\text{d}\,b^{k-1}}=e^{z\,b}\left(z^{k-1}\,H(-b)-\Gamma(k)\,\sum_{j=2}^{k}\frac{z^{k-j}}{(k-j)!}\zeta(j,-b+1)\right)
\end{equation}
\indent After some proper transformations, one concludes the result matches the formula, implying an equivalence between the differentiation and the recurrence methods.

\subsection{A new approach for the polylog} \label{polylog}

There are a few different ways to obtain the limit of \eqrefe{Poly_parc} as $n$ approaches infinity than the one created in \citena{Lerch_pos}.\\

The first one is to let $b$ go to infinity in formula \eqrefe{Lerch_full_new}, assuming $\Re(z)<0$. Before doing this, however, one must eliminate $H(b)$ from the second summation, since it diverges. This is achieved by taking the formula for the particular case $k=1$, multiplying it by $z^{k-1}/(k-1)!$ and subtracting the result from \eqrefe{Lerch_full_new}. This creates a formula that also holds for any complex $z$, except if both $\Re(z) \ge 0$ and $\abs{\Im(z)}>2\pi$:
\begin{multline} \label{eq:polylog_full_trans}
\mathrm{Li}_{k}(e^z)=-\frac{z^k}{2\,k!}+\frac{z^k}{2(k-1)!}+\frac{z^{k-1}}{(k-1)!}\mathrm{Li}_{1}(e^z)+\sum_{j=2}^{k}\frac{z^{k-j}}{(k-j)!}\zeta(j)\\
-\frac{z^{k}}{2(k-1)!}\int_{0}^{1}\left((1-u)^{k-1}-1\right)\coth{\frac{z u}{2}}\,du \text{}
\end{multline}
\indent Another way to obtain the polylog is to eliminate $H(n)$ from equation \eqrefe{Poly_parc} in the same way that was just described. One then applies the delta operator (with respect to $n$) to both sides of the resulting equation. It can then be summed over the natural numbers since no terms diverge, and this leads to an alternative expression for the polylog.

\subsection{A proof of the polylog formula}

Although the polylog formula stems from the Lerch $\Phi$ formula, there is another clever way to prove it.\\

The equation \eqrefe{polylog_full_trans}, which is proved, is a transformation of the polylog formula introduced in reference \citena{Lerch_pos},
\begin{multline} \label{eq:Poly_full}
\mathrm{Li}_{k}(e^z)=-\frac{z^{k}}{2\,k!}-\frac{z^{k-1}}{(k-1)!}\log{\left(-\frac{z}{2\pi}\right)}+\sum_{j=2}^{k}\frac{z^{k-j}}{(k-j)!}\zeta(j)\\ -\frac{z^{k}}{2(k-1)!}\int_{0}^{1}(1-u)^{k-1}\coth{\frac{z u}{2}}-\frac{2\pi}{z}(1-u)\cot{\pi u}\,du \text{,}
\end{multline}
\noindent since \eqrefe{Poly_full} can be transformed into \eqrefe{polylog_full_trans} in the same way described previously. Note that \eqrefe{Poly_full} has a slight advantage, as its domain includes $z=2\pi\ii$, whereas the domain of \eqrefe{polylog_full_trans} does not (one of the terms is not defined).\\

The two formulae, therefore, imply the following equality,
\begin{equation} \nonumber
\int_{0}^{1}z\coth{\frac{z u}{2}}-2\pi(1-u)\cot{\pi u}\,du=-z-2\log{\left(-\frac{z}{2\pi}\right)}-2\,\mathrm{Li}_{1}(e^z) \text{,}
\end{equation}
\noindent which can be evaluated more easily by differentiating both sides with respect to $z$. The differentiated integral has a closed form (obtained via software) that matches the derivative of the function on the right-hand side. Since the original equation holds at least at a point $z$, the two formulae are equivalent (where both are defined) and thus this constitutes a proof of formula \eqrefe{Poly_full}.

\section{A self-contained Lerch $\Phi$ formula}

The objective in this section is to create a stand-alone formula for the Lerch $\Phi$, that is, a formula that does not depend on the Hurwitz zeta function. To achieve this, a new easier method to generalize formulae is introduced, which consists in starting from a simple formula and then generalize it through successive differentiation with respect to the variable $b$. As will be seen, this new formula is the analytic continuation of $\Phi$ and, like in reference \citena{RFE_Ext}, a new term must be added at the end to correct for negative shift ($b$).

\subsection{The $H_k(b)$ formula}

First, one needs a formula for $H_k(b)$. For this, the method is easier to apply when equation \eqrefe{Lerch_full_new} is employed, with $z=2\pi\ii$. Then one has,
\begin{multline} \nonumber
\left(e^{2\pi\ii\,b}-1\right)H(b)=-\pi\ii+\left(e^{2\pi\ii\,b}-1\right)\left(-\frac{\pi\ii}{2}+\frac{1}{2\,b}\right)\\+\pi\,\int_0^1\left(-\left(e^{2\pi\ii\,b}-1\right)u+e^{2\pi\ii\,b\,u}-1\right)\cot{\pi u}\,du \text{,}
\end{multline}
\noindent where equation \eqrefe{Poly_full} was used to replace the singularity, $\mathrm{Li}_{1}(1)$. Therefore:
\begin{equation} \label{eq:H(b)_init}
H(b)=-\frac{\pi\ii}{2}+\frac{1}{2\,b}-\frac{\pi\ii}{e^{2\pi\ii\,b}-1}+\pi\,\int_0^1\left(-u+\frac{e^{2\pi\ii\,b\,u}-1}{e^{2\pi\ii\,b}-1}\right)\cot{\pi u}\,du \text{}
\end{equation}
\indent To obtain the generalization, each side of the equation is differentiated $k-1$ times with respect to $b$, noting that:
\begin{equation} \nonumber
\frac{e^{2\pi\ii\,b\,u}}{e^{2\pi\ii\,b}-1}=-\frac{e^{2\pi\ii\,b\,u}}{2}\left(1+\ii\cot{\pi b} \right) \text{}
\end{equation}
\indent Section \secrefe{Leibniz} provides a more detailed description of how the general derivative of the product of two functions is produced with the Leibniz rule. Particularly, the derivatives of the cotangent are given in formula \eqrefe{Cot_derivatives}. The derivative below is obtained from equation \eqrefe{e^z_cot_deriv_2}, plus the chain rule:
\begin{equation} \label{eq:e^z_cot_deriv_3}
\frac{\text{d}^{k-1}}{\text{d}\,b^{k-1}}\left(\frac{e^{2\pi\ii\,b\,u}}{e^{2\pi\ii\,b}-1}\right)=(-2\pi\ii)^{k-1}\,e^{-2\pi\ii\,b(1-u)}\,\Phi\left(e^{-2\pi\ii\,b},-k+1,1-u\right) 
\end{equation}
\indent When all are combined\footnote{The Kronecker delta is being used for the special case $k=1$.}, the formula must hold for $\Re(k)>1$ or $k=1$:
\begin{multline} \label{eq:H_k(b)}
H_k(b)=-\frac{\pi\ii}{2}\,\delta_{1,k}+\frac{1}{2\,b^k}-\ii\,\Im{\left(\frac{1}{b^k}+2\sum_{j=1}^{\floor{-b}}\frac{1}{(j+b)^k} \right)}+(1-\delta_{1,k})\,\zeta(k)-\frac{(2\pi\ii)^{k}}{2\,\Gamma(k)}\mathrm{Li}_{-k+1}\left(e^{-2\pi\ii\,b}\right)\\+\frac{\ii\,(2\pi\ii)^{k}}{2\,\Gamma(k)}\int_0^1\Big(\delta_{1,k}\,u-e^{-2\pi\ii\,b(1-u)}\,\Phi\left(e^{-2\pi\ii\,b},-k+1,1-u\right)\\+\mathrm{Li}_{-k+1}\left(e^{-2\pi\ii\,b}\right)\Big)\cot{\pi u}\,du 
\end{multline}
\indent As was the case for $\zeta(k,b)$ in reference \citena{RFE_Ext}, the formula fails when the shift is negative (and $k$ is not a positive integer), but by a known quantity, the new term that was added -- though it only works if $k$ is real, since it ensures that the imaginary part comes from the first few terms of the series.

\subsection{The summation of $H_j(b)$}

The next step is to produce the analytic continuation of the summation of $H_j(b)$ over $j$ from \eqrefe{Lerch_full_new}, which is actually equivalent to a combination of terms from \eqrefe{Lerch_full_alt}. Note that here the imaginary part is not included, since it can be deduced at the end, from the final Lerch $\Phi$ formula.\\ 

To simplify the calculation, noting that,
\begin{equation} \nonumber
\frac{1}{2\,b^{k}}\sum_{j=0}^{k}\frac{(z\,b)^j}{j!}=\frac{z^k}{2\,k!}+\frac{1}{2}\sum_{j=1}^{k}\frac{z^{k-j}}{(k-j)!}\frac{1}{b^j}\text{,}
\end{equation}
\noindent two terms of \eqrefe{Lerch_full_new} are taken together:
\begin{multline} \nonumber
\sum_{j=1}^{k}\frac{z^{k-j}}{(k-j)!}\left(-\frac{1}{2\,b^j}+H_j(b)\right)\\=\frac{z^{k-1}}{(k-1)!}\left(-\frac{1}{2\,b}+H(b)\right)+\sum_{j=2}^{k}\frac{z^{k-j}}{(k-j)!}\left(-\frac{1}{2\,b^j}+\zeta(j)-\zeta(j,b+1)\right)\\=-\frac{\pi\ii}{2 \,\Gamma(k)}+\sum_{j=2}^{k}\frac{z^{k-j}}{(k-j)!}\zeta(j)-\frac{1}{2}\sum_{j=1}^{k}\frac{z^{k-j}}{(k-j)!}\frac{(2\pi\ii)^{j}}{(j-1)!}\,\mathrm{Li}_{-j+1}(e^{-2\pi\ii\,b})\\+\frac{\ii}{2}\int_0^1 \Big(\frac{2\pi\ii\,z^{k-1}\,u}{\Gamma(k)}-\sum_{j=1}^{k}\frac{z^{k-j}}{(k-j)!}\frac{(2\pi\ii)^{j}}{(j-1)!}\Big(e^{-2\pi\ii\,b(1-u)}\,\Phi\left(e^{-2\pi\ii\,b},-j+1,1-u\right)\\-\mathrm{Li}_{-k+1}\left(e^{-2\pi\ii\,b}\right)\Big)\Big)\cot{\pi u}\,du \text{,}
\end{multline}
\noindent where a relation between $H_j(b)$ and $\zeta(j,b)$, \eqrefe{Hurwitz_zeta_rel}, was employed. An analytic continuation for the summation within the integral can be produced with the identity below, from reference \citena{Lerch_neg}:
\begin{multline} \label{eq:Soma_Lerch_neg}
\frac{1}{\Gamma(k)}\,\sum_{j=1}^{k}\binom{k-1}{j-1}\,(2\pi\ii)^{j}\,z^{k-j}\,\Phi\left(e^{-2\pi\ii\,b},-j+1,1-u\right)\\=\frac{(2\pi\ii)^{k}}{\Gamma(k)}\,\Phi\left(e^{-2\pi\ii\,b},-k+1,1-u+\frac{z}{2\pi\ii}\right) \text{}
\end{multline}
\indent The second summation, on the other hand, is a particular case of the above -- it is obtained just by multiplying by $e^{-2\pi\ii\,b(1-u)}$ and setting $u=0$.\\

When all are combined, one arrives at this formula that holds even for purely imaginary $z$ with $\Im(z)>2\pi$:
\begin{multline} \label{eq:sum_Hurwitz}
-\frac{z^{k-1}}{(k-1)!}\left(-\frac{1}{2\,b}+H(b)\right)+\sum_{j=2}^{k}\frac{z^{k-j}}{(k-j)!}\left(-\frac{1}{2\,b^j}+\zeta(j,b)\right)\\= -\frac{\pi\ii\,z^{k-1}}{2\,\Gamma(k)}+\frac{(2\pi\ii)^{k}}{2\,\Gamma(k)}\,\Phi\left(e^{-2\pi\ii\,b},-k+1,\frac{z}{2\pi\ii}\right)\\+\frac{\ii\,(2\pi\ii)^{k}}{2\,\Gamma(k)}\int_0^1 \Big(\left(\frac{z}{2\pi\ii}\right)^{k-1}(1-u)+e^{-2\pi\ii\,b(1-u)}\,\Phi\left(e^{-2\pi\ii\,b},-k+1,1-u+\frac{z}{2\pi\ii}\right)\\-\Phi\left(e^{-2\pi\ii\,b},-k+1,\frac{z}{2\pi\ii}\right)\Big)\cot{\pi u}\,du \text{}
\end{multline}

\subsection{The Lerch $\Phi$ analytic continuation}

When formula \eqrefe{sum_Hurwitz} is replaced into \eqrefe{Lerch_full_alt}, an analytic continuation of the Lerch $\Phi$ formula that holds when $\Re(k)>1$ is obtained:
\begin{multline} \label{eq:Lerch_AC}
e^{z\,b}\,\Phi(e^z,k,b)=\frac{e^{z\,b}}{2\,b^{k}}-\frac{\pi\ii\,z^{k-1}}{2\,\Gamma(k)}+\ii\,e^{z\,b}\,\Im{\left(-\frac{1}{b^k}+2\sum_{j=0}^{\floor{-b}}\frac{e^{z\,j}}{(j+b)^k} \right)}\\-\frac{z^{k-1}}{\Gamma(k)}\log{\left(-\frac{z}{2\pi}\right)}+\frac{(2\pi\ii)^{k}}{2\,\Gamma(k)}\,\Phi\left(e^{-2\pi\ii\,b},-k+1,\frac{z}{2\pi\ii}\right)\\-\frac{1}{2\,\Gamma(k)}\int_0^1 z^k\,(1-u)^{k-1}\,e^{z\,b\,u}\coth{\frac{z u}{2}}\\-\ii\,(2\pi\ii)^{k}\Big(e^{-2\pi\ii\,b(1-u)}\,\Phi\left(e^{-2\pi\ii\,b},-k+1,1-u+\frac{z}{2\pi\ii}\right)\\-\Phi\left(e^{-2\pi\ii\,b},-k+1,\frac{z}{2\pi\ii}\right)\Big)\cot{\pi u}\,du \text{,}
\end{multline}
\noindent where the extra term (the imaginary part) was found by analogy with the case presented in \citena{RFE_Ext}. That term is the correction of the formula for negative shift $b$, assuming $z$ real, and has been tested with software.\\

This formula is useful for the generalization of the Riemann functional equation that is developed later in this paper. Since the imaginary term does not appear in the initial equation from which it is derived, this formula identifies what is missing in that generalization.

\section{Generalization via recurrence} \label{AC_summation}

We now know that the summations that appear in the formulae of the Lerch $\Phi$ and the polylog functions admit analytic continuations, but can one produce closed forms for them? If the answer were yes, that would imply a closed form for the integral of \eqrefe{Lerch_AC}.\\

Having a closed form for the analytic continuation of these summations would be remarkable, since it would provide a shortcut to the derivation of the generalized functional equation of the next section.

\subsection{Polylog}

Starting from the polylog, even though it is the simplest of the cases, obtaining a closed form for it has proved very challenging.\\

The analytic continuation of the summation may be the following infinite series, 
\begin{equation} \label{eq:sum_of_zetas}
f(z)=\sum_{j=2}^{k}\frac{z^{k-j}}{(k-j)!}\zeta(j)=\frac{1}{\Gamma(k)}\sum_{j=1}^{\infty}\left(-\frac{z^k}{k}-\frac{z^{k-1}}{j}+\frac{e^{z\,j}\,\Gamma(k+1,z\,j)}{k\,j^k} \right) \text{,}
\end{equation}
\noindent where the equivalence with the upper incomplete gamma appears in reference \citena{Abra}, after the definition of the zeta function is applied and the order of the summations is inverted. However, even assuming that this is the true analytic continuation being sought, this new series is not very tractable, as it has terms that diverge.\\

Surprisingly, a better answer lies in the summation itself. The most obvious choice would be to sum the terms to infinity, but this does not work since the resulting series diverges (e.g., for a non-integer positive $k$, the terms become the reciprocal of the Gamma function at large negative arguments, and this grows faster than the decay of the powers of $z$). The solution is to express zeta at the positive integers in terms of their values at the negative integers (since they are related by the Riemann functional equation),
\begin{equation} \nonumber
\zeta(j)=\left\{
\begin{aligned}
&\frac{(2\pi\ii)^j}{2\,\Gamma(j)}\zeta(1-j) \text{, if $j$ is even}, \\
&\frac{2\,(2\pi)^{j-1}}{\Gamma(j)}\sin{\frac{\pi j}{2}}\,\zeta'(1-j) \text{, if $j$ is odd,}
\end{aligned}
\right.
\end{equation}
\noindent since in this new form it becomes possible to obtain a closed form for the summation, at least in the first case.\\

This solution exposes a divide in how one can express the function $f(z)$: it necessarily involves a separation of the summation of zeta at the even integers (for which there is a known closed from) from zeta at the odd integers (for which there is no known closed form). Nonetheless, one now has:
\begin{multline} \label{eq:sum_of_zetas2}
f(z)=\sum_{j=2}^{k}\frac{z^{k-j}}{(k-j)!}\zeta(j)=\frac{(2\pi\ii)^{k}}{2\,\Gamma(k)}\,\sum_{j=2}^{k}\binom{k-1}{j-1}\left(\frac{z}{2\pi\ii}\right)^{k-j}\zeta(1-j)\\+\frac{2\,(2\pi)^{k-1}}{\Gamma(k)}\,\sum_{j=2}^{k}\binom{k-1}{j-1}\left(\frac{z}{2\pi}\right)^{k-j}\sin{\frac{\pi j}{2}}\,\zeta'(1-j) 
\end{multline}
\indent The first summation has a closed form (a negative Hurwitz zeta), while the second one does not. The best one can do is to produce a not so convoluted power series for it, which is done in section \secrefe{Invariance}.

\section{Generalization via differentiation} \label{polylog_2_lerch}

The generating function of the polylog formula from \eqrefe{Poly_full} allows a new formula for $\Phi(e^z,1,-b+1)$ to be derived:
\begin{multline} \nonumber
\frac{1}{b}\sum_{j=1}^{\infty}b^j\,\mathrm{Li}_{j}(e^z)=-\frac{1}{2\,b}\sum_{j=1}^{\infty}\frac{(z\,b)^{j}}{j!}-\sum_{j=1}^{\infty}\frac{(z\,b)^{j-1}}{(j-1)!}\log{\left(-\frac{z}{2\pi}\right)}+\frac{1}{b}\sum_{j=1}^{\infty}b^j\sum_{j=2}^{j}\frac{z^{j-j}}{(j-j)!}\zeta(j)\\ -\frac{z}{2}\int_{0}^{1}\sum_{j=1}^{\infty}\frac{(z\,b(1-u))^{j-1}}{(j-1)!}\coth{\frac{z u}{2}}-\frac{2\pi}{z}\sum_{j=1}^{\infty}\frac{(z\,b)^{j-1}}{(j-1)!}(1-u)\cot{\pi u}\,du \text{}
\end{multline}
\indent The advantage of this formula is that it can be used to uncover a relation between the Lerch $\Phi$ at the positive and at the negative order parameters, just like the Riemann functional equation. (However, it should be noted that using $\Phi(e^z,1,b)$ from formula \eqrefe{Lerch_full} would work just the same for this purpose -- the result does not change).  The  equation above simplifies to:
\begin{multline} \label{eq:New_Lerch}
e^{z}\,\Phi(e^z,1,-b+1)=\frac{1-e^{z\,b}}{2\,b}-e^{z\,b}\log{\left(-\frac{z}{2\pi}\right)}+\frac{e^{z\,b}}{b}\,\sum_{j=2}^{\infty}b^j\,\zeta(j)\\ -\frac{z}{2}\int_{0}^{1}e^{z\,b(1-u)}\coth{\frac{z u}{2}}-\frac{2\pi\,e^{z\,b}}{z}(1-u)\cot{\pi u}\,du \text{}
\end{multline}
\indent The justification for the above is the transformation below, where the expression on the right-hand side only holds for $\abs{b}<1$:
\begin{equation} \nonumber
\sum_{k=2}^{\infty}b^k\sum_{j=2}^{k}\frac{z^{k-j}}{(k-j)!}\zeta(j)=e^{z\,b}\,\sum_{j=2}^{\infty}b^j\,\zeta(j) \text{}
\end{equation}

\subsection{Exploiting the symmetries} \label{symmetries}

If the signs of both variables in \eqrefe{New_Lerch} are reversed, the equation becomes: 
\begin{multline} \nonumber
e^{-z}\,\Phi(e^{-z},1,b+1)=-\frac{1-e^{z\,b}}{2\,b}-e^{z\,b}\log{\left(\frac{z}{2\pi}\right)}-\frac{e^{z\,b}}{b}\,\sum_{j=2}^{\infty}(-b)^j\,\zeta(j)\\ -\frac{z}{2}\int_{0}^{1}e^{z\,b(1-u)}\coth{\frac{z u}{2}}-\frac{2\pi\,e^{z\,b}}{z}(1-u)\cot{\pi u}\,du \text{}
\end{multline}
\indent Note the integral does not change, so if this new equation is subtracted from \eqrefe{New_Lerch} it vanishes:
\begin{multline} \nonumber
e^{z}\,\Phi(e^z,1,-b+1)-e^{-z}\,\Phi(e^{-z},1,b+1)\\=\frac{1-e^{z\,b}}{b}-e^{z\,b}\left(\log{\left(-\frac{z}{2\pi}\right)}-\log{\left(\frac{z}{2\pi}\right)}\right)+\frac{2\,e^{z\,b}}{b}\,\sum_{j=1}^{\infty}b^{2j}\,\zeta(2j) 
\end{multline}
\indent The ordinary generating function of the Riemann zeta function at the even integers is a well-known result from the literature\citesup{Apostol}\citesup{DLMF}:
\begin{equation} \nonumber
\sum_{j=1}^{\infty}b^{2j}\,\zeta(2j)=\frac{1-\pi b\cot{\pi b}}{2}
\end{equation}
\indent When this closed form is replaced into the formula, after a few possible simplifications it yields the two-variable equation:
\begin{multline} \label{eq:Funct_eq_1}
e^{z}\,\Phi(e^z,1,-b+1)-e^{-z}\,\Phi(e^{-z},1,b+1)\\=\frac{1}{b}-e^{z\,b}\left(\log{(-z)}-\log{(z)}\right)-\pi\,e^{z\,b}\cot{\pi b}
\end{multline}
\indent Note that taking the generating function of the Lerch $\Phi$ formula from \eqrefe{Lerch_full}, instead of the polylog, leads to essentially the same equation (except for the introduction of a redundant variable into it).

\subsection{The Leibniz rule for derivatives} \label{Leibniz}

To make the order parameter of the Lerch $\Phi$ function go from one to $k$ in \eqrefe{Funct_eq_1}, both sides of the equation are differentiated $k-1$ times with respect to $b$, which can be achieved with the Leibniz rule.\\

The Leibniz rule\citesup{Abra} is a binomial formula used to find the $k$-th derivative of the product of two functions:
\begin{equation} \label{eq:Leibniz}
\frac{\text{d}^k(f(x)\cdot g(x))}{\text{d}\,x^k}=\sum_{j=0}^{k}\binom{k}{j}f^{(j)}(x)g^{(k-j)}(x) \text{}
\end{equation}
\indent Here the only term that is not trivial to differentiate $k-1$ times is $e^{z\,b}\cot{\pi b}$. In order to employ the Leibniz rule, one can rely on this identity from reference \citena{Lerch_neg}:
\begin{equation} \label{eq:Cot_derivatives}
\frac{\text{d}^k(\cot{a\,x})}{\text{d}\,x^k}=-\ii\,\delta_{0\,k}-2\,\ii(2\,\ii\,a)^k\,\mathrm{Li}_{-k}\left(e^{2\,\ii\,a\,x}\right) \text{}
\end{equation}
\indent Therefore, differentiating the only term of \eqrefe{Funct_eq_1} that is not trivial $k-1$ times with respect to $b$ gives:
\begin{equation} \nonumber
\frac{\text{d}^{k-1}(e^{z\,b}\cot{\pi b})}{\text{d}\,b^{k-1}}=e^{z\,b}\sum_{j=0}^{k-1}\binom{k-1}{j}\left(-\ii\,\delta_{0\,j}-2\,\ii(2\pi\ii)^j\,\mathrm{Li}_{-j}\left(e^{2\pi\ii\,b}\right)\right)z^{k-1-j}
\end{equation}
\indent An analytic continuation for this summation can be produced with the identity below, again from reference \citena{Lerch_neg}:
\begin{equation} \label{eq:Soma_Polylog}
\Phi(e^{z},-k,b+1)=e^{-z}\,\sum_{j=0}^{k}\binom{k}{j}\mathrm{Li}_{-j}\left(e^{z}\right)b^{k-j} \text{,}
\end{equation}
\noindent which yields,
\begin{equation} \label{eq:e^z_cot_deriv_1}
\frac{\text{d}^{k-1}(e^{z\,b}\cot{\pi b})}{\text{d}\,b^{k-1}}=-\ii\,e^{z\,b}\left(z^{k-1}+2\,(2\pi\ii)^{k-1}\,e^{2\pi\ii\,b}\,\Phi\left(e^{2\pi\ii\,b},-k+1,1+\frac{z}{2\pi\ii}\right)\right) \text{,}
\end{equation}
\noindent or since,
\begin{equation} \label{eq:ident_Lerch}
\Phi\left(e^{2\pi\ii\,b},-k+1,1+\frac{z}{2\pi\ii}\right)=e^{-2\pi\ii\,b}\left(-\left(\frac{z}{2\pi\ii}\right)^{k-1}+\Phi\left(e^{2\pi\ii\,b},-k+1,\frac{z}{2\pi\ii}\right)\right) \text{,}
\end{equation}\\
\noindent one also has the alternative formulation below:
\begin{equation} \label{eq:e^z_cot_deriv_2}
\frac{\text{d}^{k-1}(e^{z\,b}\cot{\pi b})}{\text{d}\,b^{k-1}}=\ii\,e^{z\,b}\left(z^{k-1}-2\,(2\pi\ii)^{k-1}\,\Phi\left(e^{2\pi\ii\,b},-k+1,\frac{z}{2\pi\ii}\right)\right) 
\end{equation}

\subsection{The generalized functional equation}

To conclude the reasoning, the following formulae are needed:
\begin{equation} \label{eq:Lerch_deriv_1}
\frac{\text{d}^{k-1}\,\Phi(e^z,1,-b+1)}{\text{d}\,b^{k-1}}=\Gamma(k)\,\Phi(e^z,k,-b+1) \text{,}
\end{equation}
\noindent and,
\begin{equation} \label{eq:Lerch_deriv_2}
\frac{\text{d}^{k-1}\,\Phi(e^{-z},1,b+1)}{\text{d}\,b^{k-1}}=-(-1)^{k}\,\Gamma(k)\,\Phi(e^{-z},k,b+1) \text{}
\end{equation}\\
\indent Since differentiating the remaining terms of \eqrefe{Funct_eq_1} is relatively simple, they are skipped.\\

The new generalized functional equation is therefore:
\begin{multline} \label{eq:funct_eqt_1_final}
\Gamma(k)\left(e^z\,\Phi(e^z,k,-b+1)+(-1)^{k}\,e^{-z}\,\Phi(e^{-z},k,b+1)\right)\\= -\frac{(-1)^{k}\,\Gamma(k)}{b^k}+(2\pi\ii)^{k}\,e^{z\,b}\,\Phi\left(e^{2\pi\ii\,b},-k+1,\frac{z}{2\pi\ii}\right)\\-z^{k-1}\,e^{z\,b}\left(\pi\ii+\log{(-z)}-\log{(z)}\right) \text{,}\\
\text{where } (-1)^k=\exp(k\pi\ii) \text{, if } k \notin \mathbb{Z}
\end{multline}
\indent And from a relation similar to \eqrefe{ident_Lerch},
\begin{equation} \label{eq:ident_Lerch2}
\Phi\left(e^{z},k,b+1\right)=e^{-z}\left(-\frac{1}{b^k}+\Phi\left(e^{z},k,b\right)\right) \text{,}
\end{equation}
\noindent the generalized functional equation can be simplified to:
\begin{multline} \label{eq:funct_eqt_2_final}
\Gamma(k)\left(\Phi(e^z,k,-b)+(-1)^{k}\,\Phi(e^{-z},k,b)\right)\\= \frac{\Gamma(k)}{(-b)^k}+(2\pi\ii)^{k}\,e^{z\,b}\,\Phi\left(e^{2\pi\ii\,b},-k+1,\frac{z}{2\pi\ii}\right)\\ -z^{k-1}\,e^{z\,b}\left(\pi\ii+\log{(-z)}-\log{(z)}\right) \text{,}\\
\text{where } (-1)^k=\exp(k\pi\ii) \text{, if } k \notin \mathbb{Z}
\end{multline}
\indent Before trying to specify the domain of these relations, it is necessary to delve into some of the details. Branch is a choice of a single-valued version of a multi-valued function. For any consistent choice of principal branch of the log, the expression below,
\begin{equation} \nonumber
\pi\ii+\log{(-z)}-\log{(z)}=\begin{cases} 
0 & \text{, if } z<0 \text{ or } 0 < \Im(z) \le 2\pi,\\
2\pi\ii & \text{, if } z>0 \text{ or } -2\pi \le \Im(z) <0.\\
\end{cases} 
\end{equation}
\noindent is always a pure branch (monodromy) term, it can only take the values zero or $2\pi\ii$. Since the term is needed for purely real $z>0$ (and other $z$-domains), it can not be dropped from the equation.\\

Throughout this manuscript, the principal branch of the complex logarithm defined by, 
\begin{equation}
\log z = \log |z| + \ii\arg(z) \iff z=|z|\cdot\exp(\ii\arg(z)) \text{, } \arg(z) \in (-\pi,\pi] \text{,}
\end{equation}
\noindent is used. All expressions involving logarithms are evaluated under this convention, and branch cuts are inherited from Mathematica's implementation, which was also the software used for testing. \\

It is difficult to specify the exact domain of the newly found relations, but a naive formulation seems to be:
\begin{equation} \nonumber
\left\{
\begin{aligned}
&k \notin \mathbb{Z}_{-} \cup \{0\}, \\
&-2\pi \le \Im(z) \le 2\pi \text{ and } z \neq 0, \\
&\Im(b) > 0, \\
&\text{or } b \le 1, \\
&\text{or } 0 < \Re(b) < 1 \text{ and } -2\pi \le \Im(b) \le 0, \\
&\text{or } -1 < \Re(b) < 0 \text{ and } 0 \le \Im(b) \le 2\pi.
\end{aligned}
\right.
\end{equation}
\indent The sets given above each contain at least one point for which the equation holds -- they are not supposed to be combined into a Cartesian product.\\

Though improper, non-positive $k$ holds at the limit. The relations must hold for $\Im(z)=2\pi$, if $\Re(k)>0$. For integer $k$ (even negative), the equation may hold outside of the presumed domain. These relations have been tested numerically for these assumptions to be made, but they are neither certain nor proven, given how convoluted a three-variable function can be.

\subsection{Invariance relative to the initial equation} \label{Invariance}

The generalized relation that is obtained by differentiating the initial equation a certain number of times does not depend on it. For example, when formula \eqrefe{Lerch_full_alt} is used (with both sides multiplied by $e^{-z\,b}$), starting from $k=q$ and differentiating it $k-q$ times with respect to $b$ yields the same result as starting from $k=1$. This is easier to see by looking at the case of the term $\Phi(e^z,q,b)$:
\begin{equation} \nonumber 
\frac{\partial^{\,k-q}}{\partial\,b^{k-q}}\,\Phi(e^z,q,b)=(-1)^{k-q}\,\frac{\Gamma(k)}{\Gamma(q)}\,\Phi(e^z,k,b)
\end{equation}
\indent But not all the terms of the formula have this property on their own, some must be combined.\\

This fact can be stated as follows. For integer $q$ with $1 \le q \le k$:
\begin{multline} \label{eq:recurrence} 
\frac{\partial^{\,k-q}}{\partial\,b^{k-q}}\left( e^{-z\,b}\left(\Gamma(q)\,\sum_{j=2}^{q}\frac{z^{q-j}}{(q-j)!}\left(-\frac{1}{2\,b^j}+\zeta(j,b)\right)-z^{q-1}\left(-\frac{1}{2\,b}+H(b)\right)\right) \right)
\\=(-1)^{k-q}\,e^{-z\,b}\left(\Gamma(k)\,\sum_{j=2}^{k}\frac{z^{k-j}}{(k-j)!}\left(-\frac{1}{2\,b^j}+\zeta(j,b)\right)-z^{k-1}\left(-\frac{1}{2\,b}+H(b)\right)\right)
\end{multline}
\indent This finding finally allows one to obtain the analytic continuation of the summation from section \secrefe{AC_summation}, although by means of a power series. One just has to set $q=1$ and make use of the power series in \eqrefe{alternating_zeta}.

\subsection{Extending the domain of $b$}

In the real domain, the functional equation below,
\begin{multline} \nonumber
\Gamma(k)\left(e^z\,\Phi(e^z,k,-b+1)+(-1)^{k}\,\Phi(e^{-z},k,b)\right)\\= (2\pi\ii)^{k}\,e^{z\,b}\,\Phi\left(e^{2\pi\ii\,b},-k+1,\frac{z}{2\pi\ii}\right)-z^{k-1}\,e^{z\,b}\left(\pi\ii+\log{(-z)}-\log{(z)}\right) \text{,}
\end{multline}
\noindent which is a transformation of \eqrefe{funct_eqt_2_final}, only holds when $b \le 1$.\\

But it is possible to extend the domain to $b>1$. First, through comparison with \eqrefe{Lerch_AC} and some experimentation one finds that, 
\begin{multline} \nonumber 
\Gamma(k)\left(e^z\,\Phi(e^z,k,-b+1)+(-1)^{k}\,\Phi(e^{-z},k,b)\right)= 2\ii\,e^{z}\,\Gamma(k)\,\Im{\left(\sum_{j=0}^{\floor{b-1}}\frac{e^{z\,j}}{(j-b+1)^k} \right)}\\+(2\pi\ii)^{k}\,e^{z\,b}\,\Phi\left(e^{2\pi\ii\,b},-k+1,\frac{z}{2\pi\ii}\right)-z^{k-1}\,e^{z\,b}\left(\pi\ii+\log{(-z)}-\log{(z)}\right) \text{,}
\end{multline}
\noindent which through the identities below (one of which is Euler's reflection formula),
\begin{equation} \nonumber 
\Im{\left(\frac{1}{(-b)^k}\right)}=-\frac{\sin{\pi k}}{b^k}\, \text{, if } b>0 \text{, and } \Gamma(k)\,\sin{\pi k}=\frac{\pi}{\Gamma(1-k)} \text{,}
\end{equation}
\noindent reduces to the simpler relation below, which requires real or purely imaginary $z$:
\begin{multline} \label{eq:funct_eqt_2_final_ext}
\Gamma(k)\left(e^z\,\Phi(e^z,k,-b+1)+(-1)^{k}\,\Phi(e^{-z},k,b)\right)= -\frac{2\pi\ii}{\Gamma(-k+1)}\,\sum_{j=1}^{\floor{b}}\frac{e^{z\,j}}{(-j+b)^k}\\+(2\pi\ii)^{k}\,e^{z\,b}\,\Phi\left(e^{2\pi\ii\,b},-k+1,\frac{z}{2\pi\ii}\right)-z^{k-1}\,e^{z\,b}\left(\pi\ii+\log{(-z)}-\log{(z)}\right) \text{}
\end{multline}
\indent Note that for positive integer $k$ the summation term vanishes.

\subsection{Extending the domain of $z$}

In the real domain, a transformation of \eqrefe{funct_eqt_2_final_ext} under the substitution $z \mapsto 2\pi\ii\,z$ can also be extended to $z>1$, by analogy with the extension of the Riemann functional equation created in \citena{RFE_Ext}, without having to repeat the arduous steps taken in that reference to deduce that formula.\\

First, note that by making $k \mapsto -k+1$ and $b \mapsto z$, equation \eqrefe{gen_Riem_eq} can be transformed into:
\begin{multline} \label{eq:gen_Riem_eq_transf} 
\Gamma(k)\left(\mathrm{Li}_{k}\left(e^{2\pi\ii\,z}\right)+(-1)^k\,\mathrm{Li}_{k}\left(e^{-2\pi\ii\,z}\right)\right)\\=(2\pi\ii)^{k}\left(\zeta(-k+1,z)+\sum_{j=1}^{\floor{z}}\frac{1}{(-j+z)^{-k+1}}\right) 
\end{multline}
\indent Since this is a boundary case ($b=0$), with a little experimentation it is possible to infer the adjustment needed for the general case.\\

Therefore, if $z \ge 0$ and $b$ is real then, except for singularities, one has:
\begin{multline} \label{eq:funct_eqt_3_final_ext}
\Gamma(k)\left(e^{2\pi\ii\,z}\,\Phi(e^{2\pi\ii\,z},k,-b+1)+(-1)^{k}\,e^{-2\pi\ii\,z}\,\Phi(e^{-2\pi\ii\,z},k,b+1)\right)\\ = -\frac{(-1)^{k}\,\Gamma(k)}{b^k}-\frac{2\pi\ii}{\Gamma(-k+1)}\,\sum_{j=1}^{\floor{b}}\frac{e^{2\pi\ii\,z\,j}}{(-j+b)^k}\\+(2\pi\ii)^{k}\,e^{2\pi\ii\,z\,b}\left(\Phi\left(e^{2\pi\ii\,b},-k+1,z\right)+\sum_{j=1}^{\floor{z}}\frac{e^{-2\pi\ii\,b\,j}}{(-j+z)^{-k+1}}\right) 
\end{multline}
\indent It should be noted that even though the floor function is only defined for real numbers, both \eqrefe{funct_eqt_2_final_ext},  \eqrefe{funct_eqt_3_final_ext} and \eqrefe{gen_Riem_eq} can be extended to non-real numbers. One way to do this is, for example, by obtaining an expression for $\floor{b}$, by making $k=0$ in \eqrefe{gen_Riem_eq}. It is then replaced into the summation of \eqrefe{gen_Riem_eq} -- the summation itself must also be replaced with a closed form to yield the extension of that formula to non-real $b$.

\section{Boundary cases}

One may now ask the question: ``What does the generalized functional equation become when each of the three variables is zero?" Or more precisely, ``What is the limit at zero, if the expression is not defined at zero?" Herein, to make it directly comparable to the generalized Riemann functional equation from the literature, the transformation $z \mapsto 2\pi\ii\,z$ is applied:
\begin{multline} \label{eq:funct_eqt_3_final}
\Gamma(k)\left(e^{2\pi\ii\,z}\,\Phi(e^{2\pi\ii\,z},k,-b+1)+(-1)^{k}\,e^{-2\pi\ii\,z}\,\Phi(e^{-2\pi\ii\,z},k,b+1)\right)\\= -\frac{(-1)^{k}\,\Gamma(k)}{b^k}+(2\pi\ii)^{k}\,e^{2\pi\ii\,z\,b}\,\Phi\left(e^{2\pi\ii\,b},-k+1,z\right)\\ -(2\pi\ii\,z)^{k-1}\,e^{2\pi\ii\,z\,b}\left(\pi\ii+\log{(-\ii\,z)}-\log{(\ii\,z)}\right) 
\end{multline}
\indent Note that with this transformation, the $z$-domain (and only $z$) has changed to $0 \le \Re(z) \le 1$. Also, the rightmost term is always constant (either zero, if $0 \le \Re(z) \le 1$, or $2\pi\ii$). Like the general equation, the domain of these particular cases is not entirely known.

\subsection{The case of non-positive integer $k$}

When $k$ is a non-positive integer, the left-hand side is undefined (the Gamma function diverges but the expression in brackets goes to zero).\\

In their book on complex analysis\citesup{Stein}, Stein and Shakarchi explain that this limit can be solved using the Gamma function's residue property. Let $n$ be a non-negative integer. Then, provided that $f$ is differentiable and analytic at $k=-n$ and $f(-n)=0$, then:
\begin{equation} \label{eq:limit_theorem}
\lim_{k \to -n}\Gamma(k)\,f(k)=\frac{(-1)^n}{n!}\,f'(-n)
\end{equation}
\indent Therefore, the marginal functional equation at the non-positive integers $k=-n$ is:
\begin{multline} \label{eq:funct_eqt_k_0}
(-1)^{n}\,\frac{\partial \,\Phi(e^{2\pi\ii\,z},k,-b)}{\partial k}\bigg|_{k=-n}+\frac{\partial \,\Phi(e^{-2\pi\ii\,z},k,b)}{\partial k}\bigg|_{k=-n}
\\= -b^{n}\log{(-b)}-\pi\ii\,\Phi\left(e^{-2\pi\ii\,z},-n,b\right)+\frac{n!\,e^{2\pi\ii\,z\,b}}{(2\pi\ii)^{n}}\,\Phi\left(e^{2\pi\ii\,b},n+1,z\right)\\-\frac{n!\,e^{2\pi\ii\,z\,b}}{(2\pi\ii\,z)^{n+1}}\left(\pi\ii+\log{(-\ii\,z)}-\log{(\ii\,z)}\right)
\end{multline}
\indent Although it is difficult to specify the exact domain of this equation, it seems to include $-1 < \Re(z) < 1$ and $\Re(b) \le 1$.

\subsection{The case $z=0$}

This next particular case, $z=0$, is straightforward to derive. For this one,  formula \eqrefe{e^z_cot_deriv_1}, instead of \eqrefe{e^z_cot_deriv_2}, is used. If $k$ is not a non-positive integer:
\begin{equation} \label{eq:funct_eqt_z_0}
\Gamma(k)\left(\zeta(k,-b+1)+(-1)^{k}\,\zeta(k,b+1)\right)\\=-\frac{(-1)^{k}\,\Gamma(k)}{b^k}+(2\pi\ii)^{k}\,\mathrm{Li}_{-k+1}\left(e^{2\pi\ii\,b}\right)
\end{equation}
\indent And by means of an identity analogous to \eqrefe{ident_Lerch},
\begin{equation} \label{eq:ident_Hurwitz}
\zeta(k,b+1)=-\frac{1}{b^k}+\zeta(k,b) \text{,}
\end{equation}
\noindent an equivalent expression is obtained,
\begin{equation} \nonumber
\Gamma(k)\left(\zeta(k,-b)+(-1)^{k}\,\zeta(k,b)\right)=\frac{\Gamma(k)}{(-b)^k}+ (2\pi\ii)^{k}\,\mathrm{Li}_{-k+1}\left(e^{2\pi\ii\,b}\right) \text{,}
\end{equation}
\noindent with the same aforementioned domain.\\

Specifically, this marginal equation implies,
\begin{equation} \nonumber
\lim_{k\to 1}\Gamma(k)\left(\zeta(k,-b)+(-1)^{k}\,\zeta(k,b)\right)=-\frac{1}{b}-\pi\ii-\pi\cot{\pi b} \text{,}
\end{equation}
\noindent which is the basis for this infinite series representation of the cotangent:
\begin{equation} \nonumber
\pi \cot{\pi b}=-\frac{1}{b}+\sum_{j=0}^{\infty}\left(\frac{1}{j+b}-\frac{1}{j-b}\right) \text{}
\end{equation}
\indent To find the limit of the above at the non-positive integers $k=-n$, the result from Stein and Shakarchi\citesup{Stein} is applied again, giving:
\begin{multline} \label{eq:funct_eqt_k_0}
(-1)^{n}\,\frac{\partial \,\zeta(k,-b)}{\partial k}\bigg|_{k=-n}+\frac{\partial \,\zeta(k,b)}{\partial k}\bigg|_{k=-n}
\\= -b^{n}\log{(-b)}-\pi\ii\,\zeta\left(-n,b\right)+\frac{n!}{(2\pi\ii)^{n}}\,\mathrm{Li}_{n+1}\left(e^{2\pi\ii\,b}\right)
\end{multline}

\subsection{The generalization of the Riemann equation}

The case $b=0$ is the one that proves that this new functional equation is a generalization of the Riemann functional equation. It is either trivial, if one assumes $k<0$, or challenging, otherwise. But both cases lead to the same solution, the $b$ term vanishes and the Lerch $\Phi$ converges to the Hurwitz zeta function, on the right-hand side.\\

The assumption that $k>0$ implies that the limit below must hold:
\begin{equation} \label{eq:limit_final}
\lim_{b \to 0}-\frac{(-1)^{k}\,\Gamma(k)}{b^k}+(2\pi\ii)^{k}\,e^{2\pi\ii\,z\,b}\,\Phi\left(e^{2\pi\ii\,b},-k+1,z\right) = (2\pi\ii)^{k}\,\zeta(-k+1,z)
\end{equation}\\
\indent Therefore, if $k$ is not a non-positive integer, one has:
\begin{multline} \label{eq:funct_eqt_b_0}
\Gamma(k)\left(\mathrm{Li}_{k}\left(e^{2\pi\ii\,z}\right)+(-1)^{k}\,\mathrm{Li}_{k}\left(e^{-2\pi\ii\,z}\right)\right)\\= (2\pi\ii)^{k}\,\zeta(-k+1,z)\\ -(2\pi\ii\,z)^{k-1}\,\left(\pi\ii+\log{(-\ii\,z)}-\log{(\ii\,z)}\right) \text{,} 
\end{multline}
\noindent which is just a transformation of the generalized Riemann functional equation from the literature (which holds if $0 \le \Re(z) \le 1$)\citesup{Apostol}, except for the extra term on the right-hand side. The extra term enables the equation to hold for purely imaginary $z$ with $\Im(z)<0$. And for integer $k$ when $-1 \le \Re(z) \le 0$.\\

This marginal functional equation at the non-positive integers $k=-n$ holds for $\abs{\Re(z)}<1$ (except $n=0$, though \eqrefe{funct_eqt_3_final} may have a limit when both $k$ and $b$ approach zero):
\begin{multline} \label{eq:funct_eqt_k_0}
(-1)^{n}\,\frac{\partial \,\mathrm{Li}_{k}(e^{2\pi\ii\,z})}{\partial k}\bigg|_{k=-n}+\frac{\partial \,\mathrm{Li}_{k}(e^{-2\pi\ii\,z})}{\partial k}\bigg|_{k=-n}
\\= -\pi\ii\,\mathrm{Li}_{-n}\left(e^{-2\pi\ii\,z}\right)+\frac{n!}{(2\pi\ii)^{n}}\,\zeta(n+1,z)\\-\frac{n!}{(2\pi\ii\,z)^{n+1}}\left(\pi\ii+\log{(-\ii\,z)}-\log{(\ii\,z)}\right)
\end{multline}

\section{Acknowledgments}

The author wishes to thank all the readers of this paper, especially mathematicians that took time out of their busy schedules to read it. Special thanks to the Institute of Mathematics and Statistics of the University of S\~ao Paulo (IME-USP), that enabled me to obtain an education in mathematics and to pursue research as a personal interest. Any feedback or developments inspired by this work would be greatly appreciated.

\appendix

\section{Analytic continuation approaches}

An equation is only useful if it expresses the same quantity in different ways, otherwise it just proves the formula. The following are some approaches that were attempted to produce the analytic continuation of the function $f(z)$ from section \secrefe{AC_summation}.

\subsection{Differentiation approach}

This approach consists in using the equivalence below,
\begin{equation} \nonumber
\sum_{k=2}^{\infty}x^k\,\sum_{j=2}^{k}\frac{z^{k-j}}{(k-j)!}\zeta(j)=e^{z\,x}\,\sum_{j=2}^{\infty}x^j\,\zeta(j) \text{,}
\end{equation}
\noindent which holds if $\abs{x}<1$, and then differentiating it $k$ times with respect to the parameter $x$ and taking the result at $x=0$, that is:
\begin{equation} \nonumber
\frac{1}{k!}\,\frac{\text{d}^k}{\text{d}\,x^k}\left(e^{z\,x}\,\sum_{j=2}^{\infty}x^j\,\zeta(j)\right)\bigg|_{x=0}=\sum_{j=2}^{k}\frac{z^{k-j}}{(k-j)!}\zeta(j)
\end{equation}
\indent If the signs of both $x$ and $z$ are reversed and the resulting equations are added, a closed form is obtained:
\begin{equation} \nonumber
e^{z\,x}\left(\sum_{j=2}^{\infty}x^j\,\zeta(j)+\sum_{j=2}^{\infty}(-x)^j\,\zeta(j)\right)= e^{z\,x}\left(1-\pi\,x\cot{\pi x}\right)
\end{equation}
\indent To make it easier to calculate the $k$-th derivative of the function above, the identity below, that stems from the Leibniz rule, is used:
\begin{equation} \nonumber
\,\frac{\text{d}^k}{\text{d}\,x^k}\left(x \cdot f(x)\right)=x\,f^{(k)}(x)+k\,f^{(k-1)}(x)
\end{equation}
\indent Without repeating all the steps that have been detailed in section \secrefe{polylog_2_lerch}, the final result is:
\begin{multline} \nonumber
\frac{\text{d}^{k}e^{z\,x}\left(1-\pi\,x\cot{\pi x}\right)}{\text{d}\,x^{k}}=e^{z\,x}\Big(z^{k}-\pi\ii\,x\,z^{k}-\pi\ii\,k\,z^{k-1}+x\,(2\pi\ii)^{k+1}\,\Phi\left(e^{2\pi\ii\,x},-k,\frac{z}{2\pi\ii}\right)
\\+k\,(2\pi\ii)^{k}\,\Phi\left(e^{2\pi\ii\,x},-k+1,\frac{z}{2\pi\ii}\right)\Big) \text{,}
\end{multline}
\noindent and therefore the below is obtained:
\begin{multline} \label{eq:ez_cot_deriv_4}
\sum_{j=2}^{k}\frac{z^{k-j}}{(k-j)!}\zeta(j)+(-1)^{k}\,\sum_{j=2}^{k}\frac{(-z)^{k-j}}{(k-j)!}\zeta(j)\\=\frac{1}{k!}\left(z^{k}-\pi\ii\,k\,z^{k-1}+k\,(2\pi\ii)^{k}\,\zeta\left(-k+1,\frac{z}{2\pi\ii}\right)\right) 
\end{multline}
\indent At the positive integers $k$, this equation holds even when $\Im(z)>2\pi$. Unfortunately, this analytic continuation does not hold for both $z$ and $-z$ at the same time, which could allow a system of equations to be formed to solve for each summation as a variable (if the system thus formed were not singular). A new equation is needed. It seems that it is possible to have a closed form for this particular equation only because the values of the zeta function at the odd integers cancel out. Obtaining the same for zeta at the odd integers is either impossible or just too difficult to achieve.

\subsection{Taylor series approach}

Another way to express the function $f(z)$ from equation \eqrefe{sum_of_zetas} in a useful way is to expand the integral of formula \eqrefe{polylog_full_trans} with Taylor series.\\

From the Taylor series expansion of the cotangent\citesup{Abra}\citesup{Ahlfors},
\begin{equation} \nonumber
\coth{\frac{z u}{2}}=-\frac{2\ii}{\pi}\sum_{j=0}^{\infty}\zeta(2j)\left(-\frac{z u}{2\pi\ii}\right)^{2j-1} \text{,}
\end{equation}
\noindent it follows that
\begin{equation} \nonumber
\left((1-u)^{k-1}-1\right)\coth{\frac{z u}{2}}=-\frac{2\ii}{\pi}\sum_{j=0}^{\infty}\zeta(2j)\left(-\frac{z u}{2\pi\ii}\right)^{2j-1}\left((1-u)^{k-1}-1\right) \text{,}
\end{equation}
\noindent and from
\begin{equation} \nonumber
\int_{0}^{1}\left((1-u)^{k-1}-1\right)u^{2j-1}\,du=\left\{
\begin{aligned}
&\frac{\Gamma(k)\,\Gamma(2j)}{\Gamma(2j+k)}-\frac{1}{2j} \text{, if }j \ge 1 \text{,}\\
&-H(k-1) \text{, if } j=0 \text{ and } \Re(k)>0
\end{aligned}
\right.
\end{equation}
\noindent it follows that:
\begin{multline} \nonumber
\int_{0}^{1}\left((1-u)^{k-1}-1\right)\coth{\frac{z u}{2}}\,du\\=-\frac{2\,H(k-1)}{z}-\frac{4\,\Gamma(k)}{z}\sum_{j=1}^{\infty}\zeta(2j)\left(-\frac{z}{2\pi\ii}\right)^{2j}\frac{\Gamma(2j)}{\Gamma(k+2j)}
\end{multline}
\indent Here the final goal was to obtain a closed form, but the power series proved too challenging to express as a closed form, so this approach was abandoned.\\

In any case, when all results are combined, for integer $k \ge 1$ and $\abs{\Im(z)} \le 2\pi$, one has:
\begin{multline} \label{eq:Poly_full_power_series}
\mathrm{Li}_{k}(e^z)=-\frac{z^{k-1}}{\Gamma(k)}\log{(-z)}+\sum_{\substack{j=0 \\ j\neq k-1}}^{s}\frac{z^{j}}{j!}\,\zeta(k-j)\\ +\frac{z^{k-1}}{\Gamma(k)}\,H(k-1)+2\,z^{k-1}\sum_{j=1}^{\infty}\zeta(2j)\left(\frac{z}{2\pi\ii}\right)^{2j}\frac{\Gamma(2j)}{\Gamma(k+2j)}
\end{multline}
\indent If $z=2\pi\ii$, one must have integer $k \ge 2$, since this series does not converge for $k=1$.

\section{A new polylog relation} \label{Polylog_rel}

For the polylog, a new relation that holds for all natural $k$ and all complex $z$ (except $z$ with both $\Re(z) \ge 0$ and $\abs{\Im(z)}>2\pi$) can be obtained by exploiting the symmetry with respect to $z$ since, again, the integrals that appear in the equations for both $z$ and $-z$ are the same:
\begin{multline} \nonumber
\sum_{j=2}^{k}\frac{z^{k-j}}{(k-j)!}\zeta(j)+\frac{z^k}{(-z)^k}\,\sum_{j=2}^{k}\frac{(-z)^{k-j}}{(k-j)!}\zeta(j)\\=\frac{z^{k}}{\Gamma(k+1)}+\mathrm{Li}_{k}(e^z)+\frac{z^{k-1}}{\Gamma(k)}\log{\left(-\frac{z}{2\pi}\right)}+\frac{z^k}{(-z)^k}\left(\mathrm{Li}_{k}(e^{-z})+\frac{(-z)^{k-1}}{\Gamma(k)}\log{\left(\frac{z}{2\pi}\right)} \right)
\end{multline}
\indent Since this relation deals with complex numbers, it is presented in the safest possible form, to minimize potential issues with symbolic manipulation, though simplications may apply in some cases.\\

Now let $g(z,k)$ be the function on the right-hand side. Then one might think that the following system of equations could be formed:
\begin{equation} \nonumber
\left\{
\begin{aligned}
&f(z)+\frac{z^k}{(-z)^k}\,f(-z)=g(z,k) \\
&\frac{(-z)^k}{z^k}\,f(z)+f(-z)=g(-z,k)
\end{aligned}
\right.
\end{equation}
\indent But this system is singular for all $k$. Another equation is needed for the system to be solvable. For example, in appendix \secrefe{RFE}, a way to express the polylog in terms of $\zeta(k,b)$ and $\zeta(k,-b)$ is found by adding a new equation that can be deduced from the properties of $\zeta(k,-b)$. Now that the summation is less mysterious, it may not be impossible to achieve that.

\section{Restated Riemann functional equation} \label{RFE}

In the paper \citena{RFE_Ext}, using the symmetries of the Hurwitz zeta formula with respect to the shift parameter $b$, we proved that if $k \ne 1$ is complex and $b$ is a positive real number then:
\begin{multline} \label{eq:gen_Riem_eq} 
\zeta(k,b)=-\sum_{j=1}^{\floor{b}}\frac{1}{(b-j)^k}\\-\ii\,(2\pi)^{k-1}\,\Gamma(1-k)\left(e^{k\,\ii\,\pi/2}\,\mathrm{Li}_{-k+1}\left(e^{2\pi\ii\,b}\right)-e^{-k\,\ii\,\pi/2}\,\mathrm{Li}_{-k+1}\left(e^{-2\pi\ii\,b}\right)\right) 
\end{multline}
\indent Notice this equation has singularities at the positive integers, though it still holds at the limit (except, maybe, for $k=1$).\\

Now, to invert this formula and express the polylog as a function of Hurwitz zeta functions, one can exploit the symmetries with respect to the parameter $b$.\\

The first step is to deduce an expression for $\zeta(k,-b)$. Since,  
\begin{equation} \label{eq:comb_final} 
\zeta(k,-b)=\sum_{j=0}^{\floor{b}}\frac{1}{(j-b)^k}+\zeta(k,\floor{b}-b+1) \text{,}
\end{equation}
\noindent and since the formula of the Hurwitz zeta with shift $0<\floor{b}-b+1 \le 1$ can be written through equation \eqrefe{gen_Riem_eq}, one concludes that: 
\begin{multline} \label{eq:gen_Riem_eq_neg} 
\zeta(k,-b)=\sum_{j=0}^{\floor{b}}\frac{1}{(j-b)^k}\\+\ii\,(2\pi)^{k-1}\,\Gamma(1-k)\left(e^{-k\,\ii\,\pi/2}\,\mathrm{Li}_{-k+1}\left(e^{2\pi\ii\,b}\right)-e^{k\,\ii\,\pi/2}\,\mathrm{Li}_{-k+1}\left(e^{-2\pi\ii\,b}\right)\right) \text{}
\end{multline}
\indent From these two formulae, a system of equations can be formed, where the singularity of the original equation is removed by Euler's reflection formula. If $b>0$, then for case $b$ one has:
\begin{multline} \label{eq:gen_Riem_eq_inv_pos} 
\frac{(2\pi)^{k}}{\Gamma(k)}\,\mathrm{Li}_{-k+1}\left(e^{2\pi\ii\,b}\right)=e^{\ii\,\pi\,k/2}\left(\zeta(k,b)+\sum_{j=1}^{\floor{b}}\frac{1}{(b-j)^k}\right)\\+e^{-\ii\,\pi\,k/2}\left(\zeta(k,-b)-\sum_{j=0}^{\floor{b}}\frac{1}{(j-b)^k}\right) \text{,}
\end{multline}
\noindent and for case $-b$ one has:
\begin{multline} \label{eq:gen_Riem_eq_inv_neg}  
\frac{(2\pi)^{k}}{\Gamma(k)}\,\mathrm{Li}_{-k+1}\left(e^{-2\pi\ii\,b}\right)=e^{-\ii\,\pi\,k/2}\left(\zeta(k,b)+\sum_{j=1}^{\floor{b}}\frac{1}{(b-j)^k}\right)\\+e^{\ii\,\pi\,k/2}\left(\zeta(k,-b)-\sum_{j=0}^{\floor{b}}\frac{1}{(j-b)^k}\right)
\end{multline}

\end{document}